\documentclass[preprint,times,12pt]{amsart}
\usepackage{amsmath,amssymb,amsfonts,latexsym}
\usepackage{eucal}
\usepackage{amsthm}
\usepackage{mathrsfs}

\input amssym.def
 
\input amssym.tex

\renewcommand{\SS}{\mathbb{S}}
\newcommand{\RR}{\mathbb{R}}

\begin{document}

\title[Directions Sufficient for Steiner symmetrization]{A countable set of directions is sufficient for Steiner symmetrization}
\author[G. Bianchi]{Gabriele Bianchi} 
\address{G. Bianchi, Dipartimento di Matematica, 
Universit\`{a} di Firenze, Firenze, Italy 50134}
\email{gabriele.bianchi@unifi.it}
\author[D. Klain]{Daniel A. Klain} 
\address{D. Klain, Department of Mathematical Sciences,
University of Massachusetts Lowell,
Lowell, MA 01854, USA}
\email{Daniel\_{}Klain@uml.edu}
\address{E. Lutwak, D. Yang, G. Zhang, Polytechnic Institute of NYU,
Six Metrotech Center,
New York, NY 11201, USA}
\email{elutwak@poly.edu, dyang@poly.edu, gzhang@poly.edu}
\author[E. Lutwak]{Erwin Lutwak} 
\author[D. Yang]{Deane Yang} 
%\email{dyang@poly.edu}
\author[G. Zhang]{Gaoyong Zhang}
%\email{gzhang@poly.edu}

\begin{abstract} A countable dense set of directions is sufficient for Steiner symmetrization, 
but the order of directions matters.
\end{abstract}

% \begin{keyword} convex body, Steiner symmetrization \MSC 52A20
% \end{keyword}

\maketitle

\section*{Introduction and background}

Denote $n$-dimensional Euclidean space by $\RR^n$, and let $K$ be a compact convex subset of $\RR^n$.
Given a unit vector $u$, 
view $K$ as a family of line segments parallel to $u$.  Slide these segments along $u$ so that
each is symmetrically balanced around the hyperplane $u^\perp$.  By Cavalieri's principle, the volume
of $K$ is unchanged by this rearrangement.  
The new set, called the {\em Steiner symmetrization} of $K$ in the direction of $u$, 
will be denoted by ${\rm s}_u K $.  It is not difficult to show that 
${\rm s}_u K $ is also convex, and that ${\rm s}_u K  \subseteq {\rm s}_u L $ whenever $K \subseteq L$.  
A little more work verifies
the following intuitive assertion: if you iterate Steiner symmetrizations of $K$ through
a suitable sequence of unit directions, the successive Steiner symmetrals of $K$ will
approach a Euclidean ball in the Hausdorff topology on compact (convex) subsets of $\RR^n$.
A detailed proof of this assertion can be found in any of 
\cite[p. 98]{Egg}, \cite[p. 172]{Gru-book}, or \cite[p. 313]{Webster}, 
as well as in Section 2 below. 

Questions remain surrounding the following issue: 
Given a convex body $K$, under what more specific conditions on the sequence
of directions $u_i$ does the sequence of Steiner symmetrals 
\begin{align}
{\rm s}_{u_i} \cdots {\rm s}_{u_1} K
\label{stseq}
\end{align}
converge to
a ball?  For example, Mani \cite{Mani} has shown that, given a sequence of unit directions $u_i$ chosen uniformly
at random, the sequence~(\ref{stseq}) converges to a ball almost surely; that is, with unit probability.
An explicit algorithm for rounding out a convex body with a periodic sequence of Steiner symmetrizations
is described by Eggleston \cite[p. 98]{Egg}.  

For over 150 years Steiner symmetrization has been a fundamental tool for attacking
problems regarding isoperimetry and related geometric inequalities \cite{Gard-BM,Gard2006,Steiner,Talenti}.  
Steiner symmetrization appears in
the titles of dozens of papers (see e.g.
\cite{bia-gro,BLM,burch,CCF,CF2,CF1,Falconer,GardX,KM1,KM2,Long,Mani,McNabb,Scott})
and plays a key role in recent work such as \cite{Hab-Sch,LYZ-Orl-proj,vansch1,vansch2}.
In spite of the ubiquity of Steiner symmetrization throughout geometric analysis, 
many elementary questions about this construction remain unanswered.  
The authors of a recent paper \cite{LYZ-Orl-proj} required a sequence of Steiner
symmetrizations that rounded out a given convex body, using only directions 
drawn from a restricted dense set of directions in the unit sphere.  
Is {\em any} dense set of directions sufficient?  

While 
this more subtle fact may be derived from known results 
of a more highly technical nature (such as recent work
of Van Schaftingen \cite{vansch1,vansch2}),
it is not explicitly stated in the literature. We give
a very simple proof that is a variation of known proofs of
the standard Steiner symmetrization convergence theorem 
(such as that given in \cite{Webster}). 
Along the way, we also address the related open question:
If a given countable dense set of directions 
can be used to round out a body $K$, will this set always work,
regardless of its {\em ordering} when arranged in a sequence?

More precisely,
suppose we are given a set of directions $\Omega$ that is dense in the unit sphere.  
Is it always possible to choose directions $u_i$ from this {\em restricted} set $\Omega$ so that
the sequence~(\ref{stseq}) converges to a ball?  The short answer is {\em Yes}, 
provided the directions are chosen (and ordered) carefully.  On the other hand, it turns out that 
an arbitrary
countable dense sequence of directions may {\em fail} to accomplish this; that is, the ordering 
of the directions could make a difference.
In Section 1 we show by explicit example that, for certain orderings of the directions $u_i$, the limit
of the sequence~(\ref{stseq}) may {\em fail to exist}.  Then, in Section 2, we give an elementary proof that, 
given a convex body $K$ and a suitably ordered choice of directions $u_i$ from a dense set $\Omega$, 
the sequence~(\ref{stseq}) converges to a ball.

\section{Not every choice works}
\label{diverge}

In this section we construct a dense sequence of directions in the unit circle whose corresponding
sequence of Steiner symmetrizations {\em fails to converge} on a substantial family of
convex bodies.  While the example is given in dimension 2, it can be easily generalized to
arbitrary finite dimension.  The example that follows demonstrates the need for care when
iterating Steiner symmetrization as an infinite process.

Let $\{p_1, p_2, \cdots \}$ denote the sequence of positive prime integers.
Recall that the sum
\begin{align}
\sum_{i=1}^\infty \frac{1}{p_i}
\label{primesum}
\end{align}
diverges \cite[p. 18]{Apo}. 
For $m \geq 1$, let $u_m$ denote the unit vector in $\RR^2$ having counter-clockwise angle 
$$\theta_m = \sum_{i=1}^m \frac{\sqrt{2}}{p_i}$$
with the horizontal axis, measured in radians.   Since $\theta_m \rightarrow \infty$, while each 
successive incremental angle
$\frac{\sqrt{2}}{p_m} \rightarrow 0$, 
the unit vectors $u_m$ form a countable dense subset of the unit circle.

Meanwhile, observe that
$$\prod_{i=1}^\infty \cos \left( \frac{\sqrt{2}}{p_i} \right)  \geq 
\prod_{i=1}^\infty \left( 1 - \frac{1}{p_i^2} \right).\\
$$
Applying the Euler product formula \cite[p. 246]{Hardy-Wright}, we obtain
\begin{align*}
\left( \prod_{i=1}^\infty \cos \left(\frac{\sqrt{2}}{p_i} \right) \right)^{-1}  & \leq \;
\prod_{i=1}^\infty \left( \frac{1}{1 - \frac{1}{p_i^2}} \right) \\
&  = \; \prod_{i=1}^\infty \left( 1 + \frac{1}{p_i^2} + \frac{1}{p_i^4} + \cdots \right)
 = \; \sum_{k = 1}^\infty \frac{1}{k^2} 
 = \; \frac{\pi^2}{6},
\end{align*}
so that
\begin{align}
\prod_{i=1}^\infty \cos \left(\frac{\sqrt{2}}{p_i} \right)  \geq  \frac{6}{\pi^2}.
\label{notz}
\end{align}

Let $\ell$ be a vertical line segment, centered at the origin, of length 1.  
Apply the sequence of Steiner symmetrizations ${\rm s}_{u_m}$ to $\ell$.
Each symmetrization has the effect of projecting the previous line segment
onto the line perpendicular to $u_m$, thereby multiplying the previous length
by the next incremental cosine, 
$\cos\left(\frac{\sqrt{2}}{p_m}\right)$.  
Since the limiting
value of the product~(\ref{notz}) is strictly positive (greater than 1/2, in fact),
while the angles $\theta_m$ cycle around the circle forever, the iterated Steiner
symmetrals of $\ell$ also spin in circles forever, while approaching a limiting positive length.

In particular, the sequence of line segments
$$ \ell_m = {\rm s}_{u_m} \cdots {\rm s}_{u_1} \ell $$
has no limit.

Now let $K$ be a cigar-shaped convex body of area $\varepsilon$
containing $\ell$ as an axis of symmetry.  By the monotonicity of Steiner symmetrization, 
each element in the sequence of Steiner symmetrals
$$ K_m = {\rm s}_{u_m} \cdots {\rm s}_{u_1} K $$
must contain the corresponding symmetral $\ell_m$, so
that the diameter of each $K_m$ exceeds $\frac{6}{\pi^2}$.
Since each $K_m$ has the same area $\varepsilon$ as the original body $K$, 
which could be made arbitrarily
small beforehand, it follows that the sequence $K_m$ cannot approximate a ball.
Indeed, for $\varepsilon < \frac{9}{\pi^3}$ the sequence $K_m$ has no limit, since the diameter
line revolves forever, but does not shrink enough to accomodate the tiny given area $\varepsilon$.

We have shown that a countable dense sequence of directions does not necessarily
lead to a well-defined limiting Steiner symmetral.

In this specific example we used the divergent series~(\ref{primesum}) as a starting point
for computational convenience.  
Gronchi \cite{Gronchi-Example} has shown that a more general family of examples can be constructed
starting with any decreasing sequence of incremental angles $\theta_i$ provided that 
$\sum_{i=1}^\infty \theta_i^2$ converges and 
$\sum_{i=1}^\infty \theta_i$  diverges.
Iterated Steiner symmetrization in the resulting sequence of directions,
applied to a sufficiently eccentric ellipse, results in a sequence of 
ellipses whose principal axes rotate forever without converging to a circle.

\section{There is always an order of directions that works}

In view of the previous example, it is necessary to show that, given a 
countable dense set of directions in $\RR^n$, it is indeed possible to
construct a sequence of directions from this set so that successive Steiner
symmetrals of a given convex body $K$ converge to a Euclidean ball.

Let $\Omega$ be a dense subset of the Euclidean unit sphere $\SS^{n-1}$.   If 
a convex body $K$ in $\RR^n$ has
volume zero, then $K$ lies in a proper subspace of $\RR^n$.  Steiner symmetrization
of $K$ in any direction outside the affine hull of $K$ has the same effect as
orthogonal projection.  One can choose directions close to, but outside, the affine hull
so as to shrink $K$ by any positive factor, along any direction inside the affine hull.
A suitable iteration will shrink the diameter of successive symmetrizations (projections)
to zero, so that symmetrals converge to a point (a Euclidean ball of radius zero).

For a convex body $K$ in
$\RR^n$ of positive volume, let $r_K=\max_{x\in K} |x|$, which is the minimal radius
of balls centered at the origin that contain $K$. Let $r_1$ be the
infimum of all $r_C$, where $C$ is obtained from finitely many successive
Steiner symmetrizations of $K$ in directions that belong to
$\Omega$. Then there is a sequence of such convex bodies $C_i$ so
that $r_{C_i}\to r_1$. Obviously, the sequence $\{C_i\}$ is bounded,
because each $C_i \subseteq r_KB$, where $B$ is the unit ball. By the
Blaschke selection theorem \cite{red,Webster}, there is a subsequence $C_{i_k}$ that
converges to a convex body $K_1$, where $r_{K_1} = r_1$. Denote $r_1B$
by $B_1$, so that $K_1 \subseteq B_1$.

We claim that $K_1=B_1$.
Assume it is not true. There is a small cap $U$ on $\partial B_1$ so that $U\cap K_1=\emptyset$.
For any line $\xi$ 
such that $\xi\cap U\neq\emptyset$, either $\xi \cap K_1 = \emptyset$ or 
the line $\xi$
intersects a longer chord in $B_1$ than in $K_1$; that is, $|\xi\cap B_1| > |\xi\cap K_1|$. 
After taking a Steiner symmetrization ${\rm s}_u K_1$ for some $u\in \Omega$,
the symmetral ${\rm s}_u K_1$ fails to intersect both $U$
and a new cap $U'$ given by the reflection of $U$ with respect to the hyperplane $u^\perp$.
Since $\Omega$ is dense in $\SS^{n-1}$, one can
continue to take symmetrizations with respect to an appropriate 
finite family of hyperplanes with normals 
$v_1, \ldots, v_s \in \Omega$
that generate finitely many caps covering the whole sphere $\partial B_1$ and generate a convex body
$K_2$ so that $|\xi\cap B_1| > |\xi\cap K_2|$ for any line such that $\xi\cap \partial B_1\neq\emptyset$.
Thus, $r_{K_2}<r_1$.   

Denote $\tilde{C}_{i_k} = s_{\rm v_s} \cdots s_{\rm v_1} C_{i_k}$.
Since $C_{i_k} \rightarrow K_1$, while Steiner symmetrization is continuous on convex bodies with 
nonempty interior
\cite[p. 312]{Webster}, we have 
$$\tilde{C}_{i_k} = s_{\rm v_s} \cdots s_{\rm v_1} C_{i_k} \longrightarrow 
s_{\rm v_s} \cdots s_{\rm v_1} K_1 = K_2.$$
Since $r_{\tilde{C}_{i_k}} \rightarrow r_{K_2}$, it follows from the definition of $r_1$
that $r_{K_2} \geq r_1$, a contradiction.

We have shown that for any convex body $K$ of positive volume there are $u_1, \ldots, u_{i_1} \in \Omega$
so that the Hausdorff distance $d$ between ${\rm s}_{u_{i_1}} \cdots {\rm s}_{u_1} K$ and
the centered ball $B_1$ with the same volume of $K$ can be arbitrarily small.

Denote by $d(K_1,K_2)$ the Hausdorff distance between compact convex sets $K_1,K_2 \subseteq \RR^n$.
For a sequence of positive numbers $\varepsilon_k\to 0$, there are $u_1, \cdots, u_{i_1} \in \Omega$
so that $d(D_1, B_1)<\varepsilon_1$, where $D_1={\rm s}_{u_{i_1}} \cdots {\rm s}_{u_1} K$.
Similarly, \sloppy there are $u_{i_1+1}, \cdots, u_{i_2}\in \Omega$ so that $d(D_2, B_1)<\varepsilon_2$,
where $D_2={\rm s}_{u_{i_2}} \cdots {\rm s}_{u_{i_1+1}} D_1$. In general,
there are  $u_{i_{k-1}+1}, \cdots, u_{i_k}\in \Omega$ so that $d(D_k, B_1)<\varepsilon_k$,
where $D_k={\rm s}_{u_{i_k}} \cdots {\rm s}_{u_{i_{k-1}+1}} D_{k-1}$. 
Since $d({\rm s}_u K, B_1) \le d(K,B_1)$ for any $K$ when $d(K,B_1)< r_{B_1}$,
there is
a sequence $K_i={\rm s}_{u_i} \cdots {\rm s}_{u_1} K \to B_1$, where $u_i\in \Omega$.

\section{Related open questions}

In Section~\ref{diverge} we described a convex body $K$ and a sequence
of directions $u_i$ for which the sequence of Steiner symmetrals
\begin{align*}
K_i = {\rm s}_{u_i} \cdots {\rm s}_{u_1} K
\end{align*}
failed to converge
in the Hausdorff topology. However, some of the examples described in Section~\ref{diverge} clearly converge in {\em shape}: there is a corresponding sequence of isometries $\psi_i$ such that
the sequence $\{ \psi_i K_i \}$ converges. Is this always the case?
If so, and supposing also that the sequence $\{u_i\}$ 
is dense in the unit sphere $\SS^{n-1}$, is the limit of the convergent sequence
$\{ \psi_i K_i \}$ always an ellipsoid?  Moreover, what happens if $K$ is permitted to be an arbitrary
(possibly non-convex) compact set?

% \bibliographystyle{amsplain}
% \bibliography{all}

% \providecommand{\bysame}{\leavevmode\hbox to3em{\hrulefill}\thinspace}
% \providecommand{\MR}{\relax\ifhmode\unskip\space\fi MR }
% %\MRhref is called by the amsart/book/proc definition of \MR.
% \providecommand{\MRhref}[2]{%
%   \href{http://www.ams.org/mathscinet-getitem?mr=#1}{#2}
% }
% \providecommand{\href}[2]{#2}

\bibliographystyle{elsarticle-num}

\end{document}